\documentclass{amsart}
\usepackage{graphicx}

\newtheorem{thm}{Theorem}
\newtheorem{lemma}[thm]{Lemma}

\newtheorem{conj}[thm]{Conjecture}

\newcommand{\pr}{\noindent \textbf{Proof: }} 
\newcommand{\bx}{\qed} 

\newcommand{\A}{{\mathbf A}}
\newcommand{\D}{{\mathcal D}}

\title{On the Harary-Kauffman Conjecture and Turk's Head Knots}
\author[N.E.\ Dowdall, T.W.\ Mattman, K.\ Meek, and P.R.\ Solis]{Nicholas E.\ Dowdall, Thomas W.\ Mattman, Kevin Meek, and Pablo R.\ Solis}
\address{Department of Mathematics, Sonoma State University,
Rohnert Park, CA 94928}
\email{nickdowdall@yahoo.com}
\address{Department of Mathematics and Statistics,
California State University, Chico,
Chico, CA 95929-0525}
\email{TMattman@CSUChico.edu}
\address{Department of Mathematics,
Florida State University,
208 Love Building,
Tallahasee, FL 32306-4510}
\email{xboogerx@hotmail.com}
\address{Department of Mathematics,
Massachusetts Institute of Technology,
Headquarters Office,
Building 2, Room 236,
77 Massachusetts Avenue,
Cambridge, MA 02139-4307
}
\email{psolis@mit.edu}
\thanks{The research was supported in part by NSF REU award 0354174 and by the MAA's NREUP program with funding from the NSF, NSA, and Moody's.}

\begin{document}

\begin{abstract}
The $m,n$ Turk's Head Knot, $THK(m,n)$, is an ``alternating $(m,n)$ torus knot.'' We prove the Harary-Kauffman conjecture for all $THK(m,n)$ except for the case where $m \geq 5$ is odd and $n \geq 3$ is relatively prime to $m$. We also give evidence in support of the conjecture in that case. 
Our proof rests on the observation that none of 
these knots have prime determinant except for $THK(m,2)$ when $P_m$ is a Pell prime.  
\end{abstract}

\maketitle

\section{Introduction}
We investigate the Harary-Kauffman~\cite{HK} conjecture for a class of knots that, following~\cite{NY}, we call the Turk's Head Knots.

\begin{conj}[Harary-Kauffman] 
Let $\D$ be an alternating knot diagram with no nugatory crossings. If the determinant of $\D$ is a prime number $p$, then every non-trivial Fox $p$-coloring of $\D$ assigns different colors to different arcs of $\D$.
\end{conj}

\begin{figure}[ht]
\begin{center}
\includegraphics[width=4.5in]{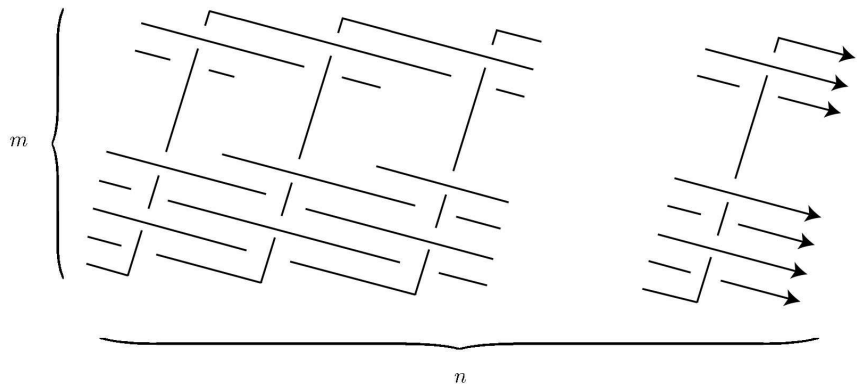}
\caption{The $THK(m,n)$ Turk's Head Link}\label{figTHKmn}
\end{center}
\end{figure} 

The link $THK(m,n)$ (where $m,n$ are both integers greater than $1$) can be formed by taking a braid representation of the
$(m,n)$ torus link and making it alternate as illustrated in figure~\ref{figTHKmn}. As is the case with torus links, this results in a link of $\mbox{GCD}(m,n)$ components; in particular, $THK(m,n)$ is a knot precisely when $m$ and $n$ are relatively prime. 

If $m = 2$ and $n$ is odd, $THK(2,n)$ is a $(2,n)$ torus knot and is also a rational knot and a Montessinos knot. 
The Harary-Kauffman Conjecture is known to hold for
such knots \cite{KL,APS}. So, we will assume that $m$ is at least $3$.  
Our key observation is that very few Turk's Head Knots have prime determinant and, therefore, most $THK(m,n)$ satisfy the Harary-Kauffman Conjecture in a trivial way. 
We conjecture that if $THK(m,n)$ is a Turk's Head Knot of prime determinant, then 
$n = 2$. For these knots the determinant is a Pell number:

\begin{thm} 
\label{thmdetm2}%
Let $m \geq 3$. The determinant of $THK(m,2)$ is the $m$th Pell number $P_m$ where $P_1 = 1$, $P_2 = 2$, and $P_m = 2P_{m-1} + P_{m-2}$ for $m \geq 3$.
\end{thm}

Thus, we conjecture that the only $THK(m,n)$ knots with prime determinant
are the $THK(m,2)$ for which $P_m$ is a Pell prime. In particular, this means that $m$ must be a prime. 
Moreover, we demonstrate that every non-trivial coloring of the diagram of $THK(m,2)$ assigns different colors to different arcs:

\begin{thm}
\label{thmKHm2}%
Let $m \geq 3$ be odd. The Harary-Kauffman conjecture holds for the diagram of $THK(m,2)$.
\end{thm}

To be precise,
theorem~\ref{thmKHm2} verifies the conjecture for the diagram of $THK(m,2)$ shown in figure~\ref{figTHKmn}. However, since minimal diagrams are related  by flypes~\cite{MT}, figure~\ref{figTHKmn} is essentially the unique minimal diagram for this knot.

Thus, in addition to the case $n=2$,
we have a proof of the Harary-Kauffman conjecture for all $THK(m,n)$ 
with composite determinant. We can prove this class includes all knots
with $m=3$ or $m$ even:

\begin{thm}
\label{thmdet}%
Let $m \geq 3$ and $n \geq 3$ be relatively prime integers.
If $m = 3$ or $m$ is even, then the Turk's Head Knot $THK(m,n)$ 
has composite determinant.
\end{thm}

For the remaining knots, that is, $THK(m,n)$ with $m \geq 5$ odd and
$n \geq 3$ relatively prime to $m$, 
we propose a formula $G_{m,n}$ for the determinant in section 6 below,
where we also show that $G_{m,n}$ is composite:

\begin{thm}
\label{thmGmn}%
$G_{m,n}$ is composite when $m,n \geq 3$ and $m$ is odd. 
\end{thm}

The structure of this paper is as follows. In the next section we introduce Fox coloring, the determinant, and connections with spanning trees of the
checkerboard graph of a knot. In sections 3, 4, 5, and 6, we prove theorems \ref{thmdetm2}, \ref{thmKHm2}, \ref{thmdet}, and \ref{thmGmn}, respectively.

\section{Coloring, determinants, and spanning trees}

In this section we give a brief overview of the idea of Fox $p$-coloring and its connections with the knot's determinant. A more complete discussion can be found in \cite{L}.

In a $p$-coloring of a knot diagram, we assign  to each arc an integer 
between $0$ and $p-1$ such that at each crossing the over strand $x$ and the two under strands $y$ and $z$ satisfy the relation $2x - y - z \equiv 0 \pmod{p}$. Further, at least two different ``colors'' must be used.

It is straightforward to show that this is a knot invariant; we will say that a knot is $p$-colorable if it has a diagram that can be 
$p$-colored. The determinant of a knot determines the valid choices for
$p$:  a knot is $p$-colorable if and only if $p$ has a common factor with the determinant of the knot. 

For an alternating knot the 
determinant is given by the number of spanning
trees of a checkerboard graph of an alternating diagram of the knot (see  \cite{BZ}). The checkerboard graph is obtained by first shading alternate
regions of the knot diagram. We then place a vertex in each shaded region
and an edge between regions that share a crossing. Figures~\ref{figspanmodd}, \ref{figspanmeven}, \ref{figwn}, and \ref{figspanmn} give examples of checkerboard graphs. See also~\cite{HK} for an introduction 
to this idea.

\section{Proof of Theorem~\ref{thmdetm2}}

\setcounter{thm}{1}

In this section we prove

\begin{thm} 
Let $m \geq 3$ be odd. The determinant of $THK(m,2)$ is the $m$th Pell number $P_m$.
\end{thm}

\setcounter{thm}{6}

\begin{figure}[ht]
\begin{center}
\includegraphics[width=2.7in]{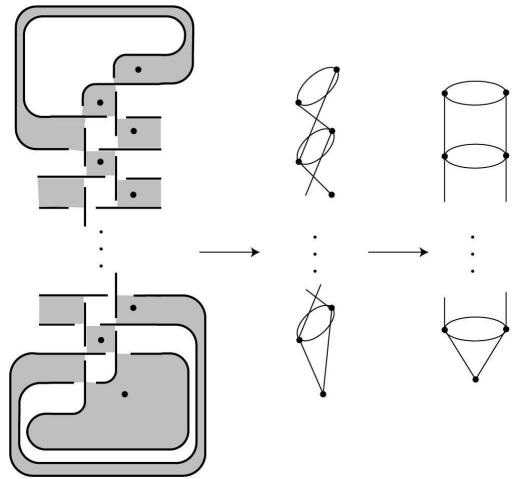}
\caption{When $m$ is odd, the graph of $THK(m,2)$ has
$m$ vertices.}\label{figspanmodd}
\end{center}
\end{figure} 

\begin{figure}[ht]
\begin{center}
\includegraphics[width=2.75in]{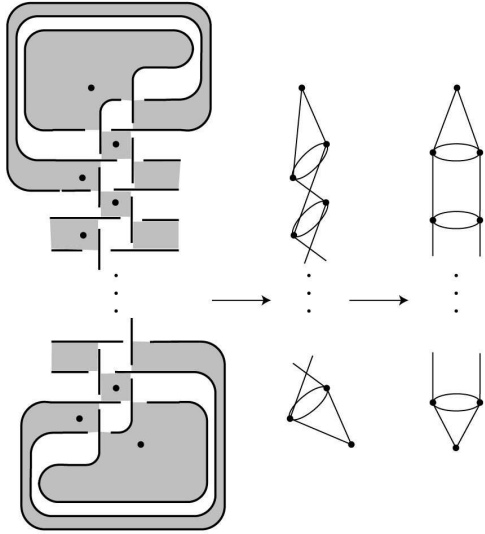}
\caption{When $m$ is even, the graph of $THK(m,2)$ again has
$m$ vertices.}\label{figspanmeven}
\end{center}
\end{figure} 

\pr We use induction to show that $P_m$ is the number of spanning trees of the checkerboard graph $H_m$ of $THK(m,2)$. The general form of the checkerboard graph $H_m$ breaks down into two cases.
If $m$ is odd we have the graph on $m$ vertices shown in figure~\ref{figspanmodd}. If $m$ is even, as shown in figure~\ref{figspanmeven},  the graph again has $m$ vertices.  First observe that the number of spanning trees of $H_{3}$ is $5=P_{3}$ and the number of spanning trees of $H_{4}$ is $12=P_{4}$.  So the theorem holds for $m=3$ and $m=4$.  Now assume that it holds for all $3 \leq m < k$.  We will show that it holds for $m = k$.

\begin{figure}[ht]
\begin{center}
\includegraphics[width=2.75in]{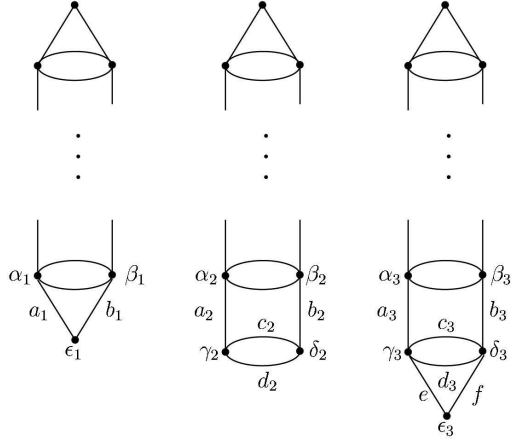}
\caption{The graphs $H_{k-2}$, $H_{k-1}$, and $H_{k}$ where $k$ is even.}\label{figm2keven}
\end{center}
\end{figure} 

First, assume that $k$ is even.  We will label vertices and edges for $H_{k-2}$, $H_{k-1}$, and $H_{k}$ as in figure~\ref{figm2keven}. Given any spanning tree of $H_{k-2}$, construct a subgraph of $H_{k}$ as follows.  All of the unlabeled edges will remain identical.  If the spanning tree for $H_{k-2}$ has edge $b_{1}$, but not $a_{1}$, we will remove $b_{1}$ and add $b_{3}$, $e$, and $f$.  If the spanning tree for $H_{k-2}$ has edge $a_{1}$, but not $b_{1}$, then we will remove $a_{1}$ and add $a_{3}$, $e$, and $f$. Finally, if the spanning tree for $H_{k-2}$ has both $a_{1}$ and
$b_{1}$, then we will remove both of these edges and add $a_{3}$, $e$, $f$, and $b_{3}$.  The newly added edges will connect $\gamma_{3}$, $\epsilon_{3}$, and $\delta_{3}$ to the rest of the graph.  All unlabeled vertices are already spanned by the original tree from $H_{k-2}$. Also, we have constructed our new graphs to avoid closed loops, so all of our new subgraphs of $H_{k}$ are spanning trees of $H_{k}$. Furthermore, in order to connect $\epsilon_{1}$ to the rest of any spanning tree of $H_{k-2}$, the spanning tree must be of one of the above forms.   Thus, for any spanning tree in $H_{k-2}$, we have exactly one associated spanning tree for $H_{k}$.  Denote by $S_{k-2}$ the set of spanning trees obtained in this manner from the set of spanning trees of $H_{k-2}$. Now, for any spanning tree of $H_{k-1}$, we can get one new subgraph of $H_{k}$ by adding the edge $e$.  We can get a second subgraph by instead adding the edge $f$.  In either case, the new subgraph connects $\epsilon_{3}$ via the new edge.  All other vertices were already connected by the spanning tree of $H_{k-1}$. Furthermore, we add no closed loops by adding exactly one of $e$ and $f$.  Thus, for every spanning tree of $H_{k-1}$, we can associate two unique spanning trees for $H_{k}$.  Denote by $S_{k-1}$ the set of spanning trees obtained in this manner.

Now we have two sets of spanning trees for $H_{k}$ with the properties that $|S_{k-2}|$ is the number of spannning trees for $H_{k-2}$ and $|S_{k-1}|$ is the twice the number of spanning trees for $H_{k-1}$.  We need to show that these two sets are mutually exclusive and that their union is the set of spanning trees of $H_{k}$.  To see that these two sets are mutually exclusive, it suffices to notice that every element of $S_{k-2}$ has $e$ and $f$ as edges, whereas every element of $S_{k-1}$ has exactly one of $e$ and $f$ as an edge.  Now if $x$ is a spanning tree of $H_{k}$, it must contain an edge which connects $\epsilon_{3}$ to the rest of the graph, so it must contain $e$ or $f$.  If $x$ contains exactly one of $e$ and $f$, then it must be in $S_{k-1}$. If $x$ contains both $e$ and $f$, then it must also contain $a_{3}$ or $b_{3}$, and must be in $S_{k-2}$.  So the number of spanning trees of $H_{k}$ is twice the number of spanning trees of $H_{k-1}$ added to the number of spanning trees of $H_{k-2}$.  But by inductive hypothesis, this means that the number of spanning trees of $H_{k}$ is $2P_{k-1}+P_{k-2}=P_{k}$, as needed.

\begin{figure}[ht]
\begin{center}
\includegraphics[width=2.8in]{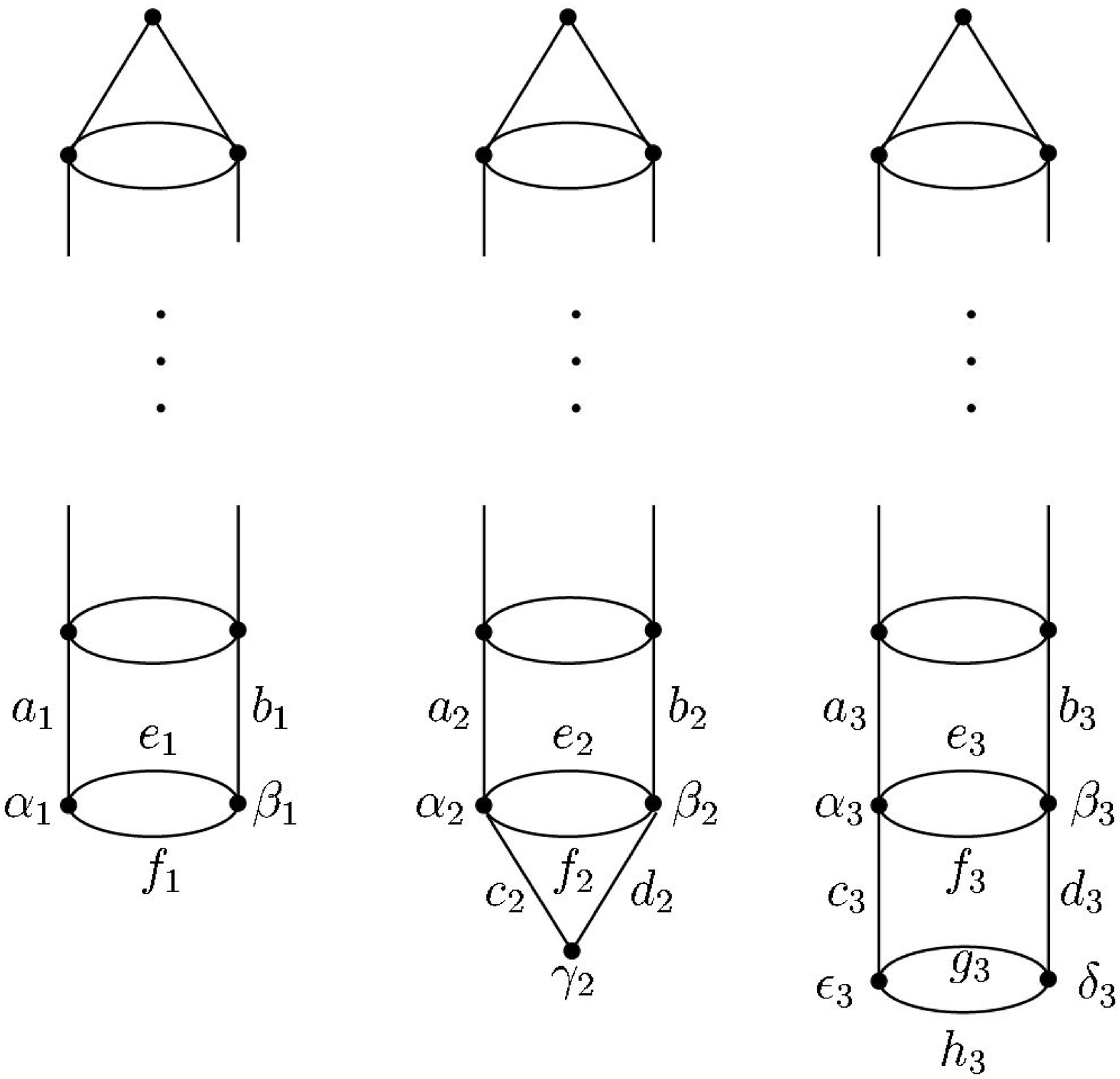}
\caption{The graphs $H_{k-2}$, $H_{k-1}$, and $H_{k}$ where $k$ is odd.}\label{figm2kodd}
\end{center}
\end{figure} 

Now we will look at the case $k$ odd.  Again, we will label graphs $H_{k-2}$, $H_{k-1}$, $H_{k}$ as in figure~\ref{figm2kodd}.  Given any spanning tree of $H_{k-2}$, we may construct a subgraph of $H_{k}$ by adding in the edges $c_{3}$ and $d_{3}$.  This will connect $\epsilon_{3}$ and $\delta_{3}$ to the rest of the graph.  We know that all other vertices are connected since we started with a spanning tree for $H_{k-2}$.  There is no way for these two new edges to create a closed loop in our subgraph, so for every spanning tree in $H_{k-2}$, we can assign a unique spanning tree of $H_{k}$.  We will denote the set of such spanning trees of $H_{k}$ by $S_{k-2}$.

To every spanning tree of $H_{k-1}$, we can associate two distinct subgraphs of $H_{k}$ as follows.  If our spanning tree contains $c_{2}$, but not $d_{2}$, we remove this edge and add $c_{3}$ and either $g_{3}$ or $h_{3}$.  If our spanning tree contains $d_{2}$ but not $c_{2}$, we replace $d_{2}$ by $d_3$ and add either $g_{3}$ or $h_{3}$.  If our spanning tree contains both $c_{2}$ and $d_{2}$, then we remove these edges and add $c_{3}$, $d_{3}$, and either $g_{3}$ or $h_{3}$.  Notice that, in the last case, we will not create a closed loop since this would imply that $c_{2}$ and $d_{2}$ would create a closed loop in our original spanning tree. Furthermore, each of these cases will connect $\epsilon_{3}$ and $\delta_{3}$ to the rest of the graph.  Finally, note that every spanning tree for $H_{k-1}$ must contain $c_{2}$ or $d_{2}$ in order to connect $\gamma_2$ to the rest of the graph.  So every spanning tree of $H_{k-1}$ can be associated with two spanning trees of $H_{k}$.  We will call the set of spanning trees of $H_{k}$ so obtained $S_{k-1}$.

Now we have two sets of spanning trees for $H_{k}$.  Furthermore, $|S_{k-2}|$ is the number of spanning trees of $H_{k-2}$ and $|S_{k-1}|$ is twice the number of spanning trees of $H_{k-1}$. We must show that these two sets are mutually exclusive, and that their union is the set of all spanning trees for $H_{k}$. For a subgraph of $H_{k}$ to be a spanning tree, $\epsilon_{3}$ and $\delta_{3}$ must be connected to the rest of the graph, so any given spanning tree must contain both $c_{3}$ and $d_{3}$, but not $g_{3}$ or $h_{3}$, placing it in $S_{k-2}$; both $c_{3}$ and $d_{3}$ as well as either $g_{3}$ or $h_{3}$ but not both, placing it in $S_{k-1}$; or it must contain exactly one of $c_{3}$ and $d_{3}$ and exactly one of $g_{3}$ and $h_{3}$, placing it in $S_{k-1}$. So, every spanning tree of $H_{k}$ is in one of these two sets.  Furthermore, in order for a graph to be in $S_{k-2}$, it cannot contain $g_{3}$ or $h_{3}$, yet every graph in $S_{k-1}$ contains exactly one of these edges.  Thus, these two sets are mutually exclusive.  In other words, the number of spanning trees of $H_{k}$ is $|S_{k-1}|+|S_{k-2}|=2P_{k-1}+P_{k-2}=P_{k}$ as needed. \bx

\section{Proof of Theorem~\ref{thmKHm2}}

\setcounter{thm}{2}

In this section we  prove

\begin{thm}
Let $m \geq 3$ be odd. The Harary-Kauffman conjecture holds for the diagram of $THK(m,2)$.
\end{thm}

\setcounter{thm}{5}

Note that we require that $m$ be odd so that $THK(m,2)$ is a knot rather than a link of two components.

\begin{figure}[ht]
\begin{center}
\includegraphics[width=2.6in]{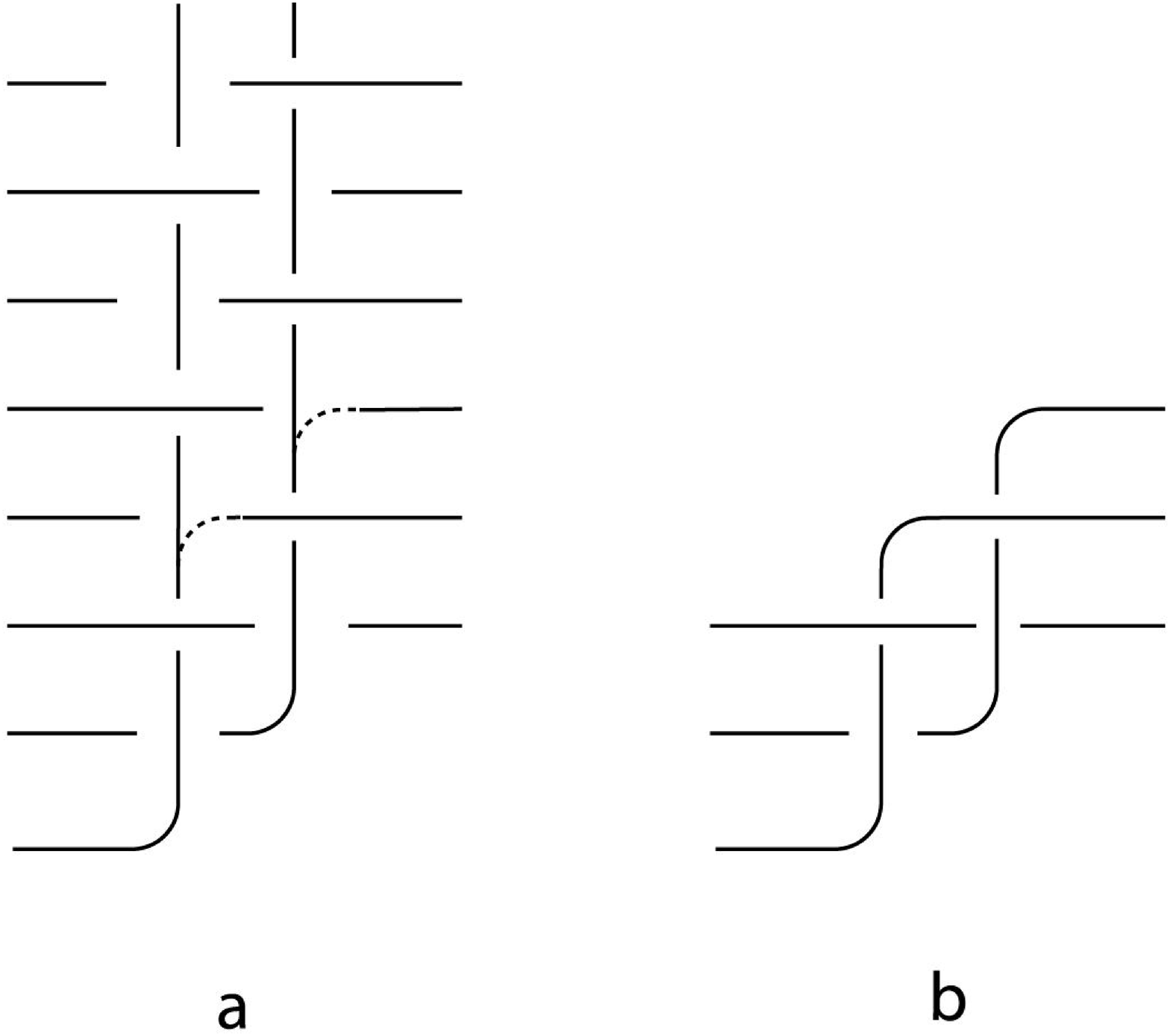}
\caption{a) The braid $THK(m,2)$. b) Connecting the dashed lines yields $THK(3,2)$}\label{figm2}
\end{center}
\end{figure} 

We will work with the general braid projection in figure~\ref{figm2}a.  To extract $THK(3,2)$, for example,  from this diagram, we need only connect two of the vertical strands with horizontal strands using the dashed lines shown.  If we have labeled these strands in a $p$-coloring this gives rise to a constraint that the colors of the newly connected strands agree modulo $p$.

We break the proof of theorem~\ref{thmKHm2} into a series of lemmas.
Our first observation is that there is essentially only one non-trivial way to color a $THK(m,2)$ diagram.
Without loss of generality, color the bottom strand at left $0$ and 
the strand above it $a$. As can be seen in figure~\ref{figacol}, 
this will uniquely determine the coloring of all
strands of the $THK(m,2)$ diagram up to the choice of $a$.
Moreover,
the $a$-coefficients come in pairs that differ in sign and 
differ by one in absolute value:

\begin{figure}[ht]
\begin{center}
\includegraphics[width=2.5in]{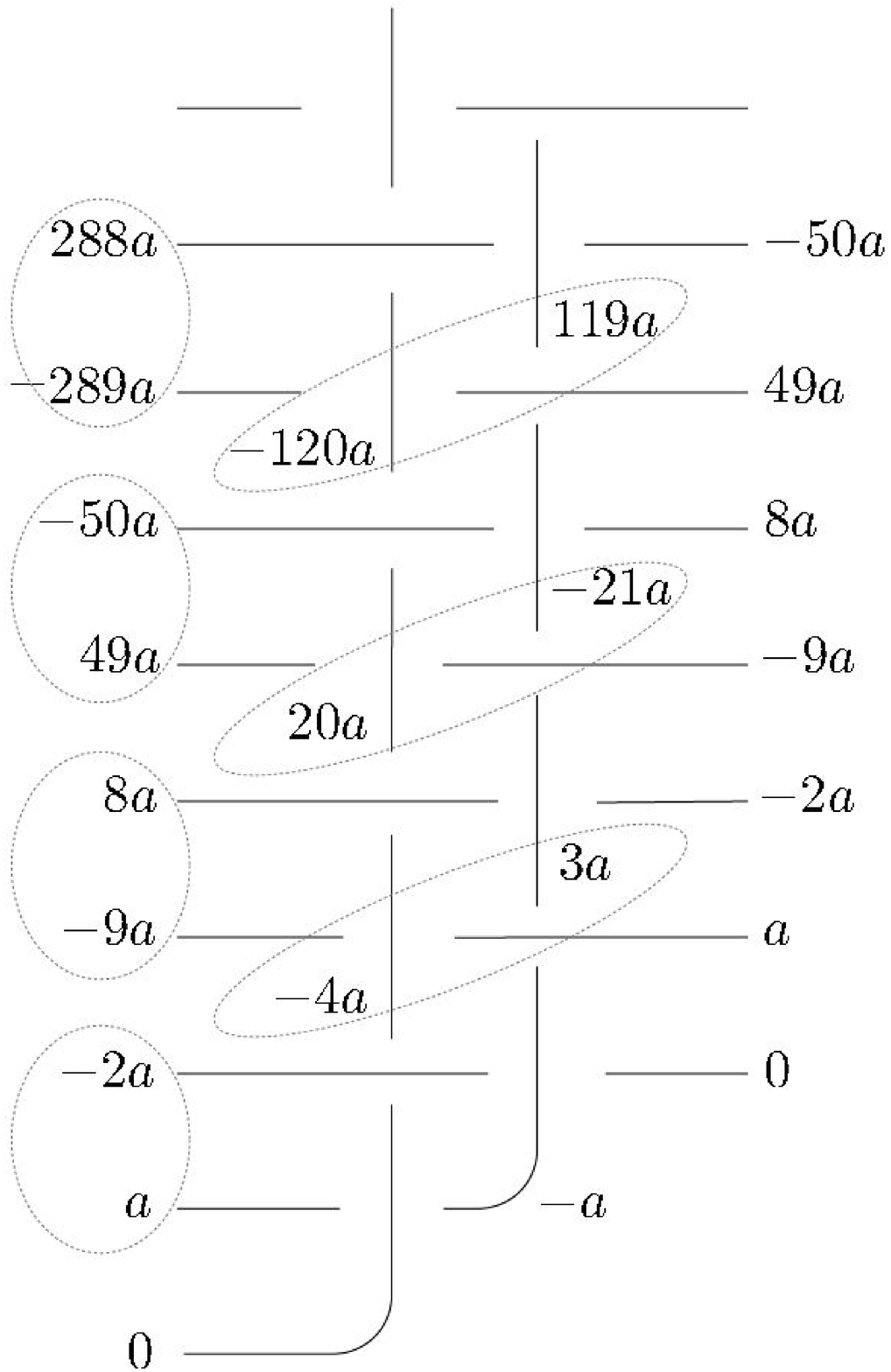}
\caption{Coloring of $THK(m,2)$ diagram.}\label{figacol}
\end{center}
\end{figure} 

\begin{lemma} \label{lemSn}

The $a$-coefficients of entering strands form adjacent pairs of opposite sign such that they differ by one in absolute value.  Coefficients of interior strands also form pairs of opposite sign such that the absolute value of the negative strand is one more than the absolute value of the positive strand. Furthermore, if we order the pairs as in figure~\ref{figacol} and choose the positive element $S_n$ from each pair, we get the following recursive relationship: $S_{n}=2S_{n-1}+S_{n-2}+1$.

\end{lemma}

\begin{figure}[ht]
\begin{center}
\includegraphics[width=3.5in]{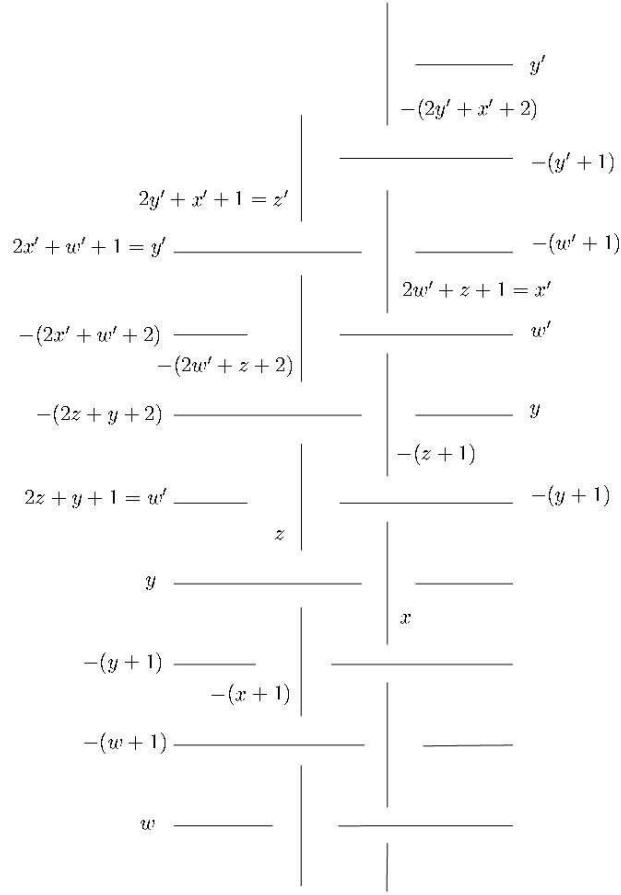}
\caption{The induction step.}\label{figkind}
\end{center}
\end{figure}

\pr 
Specifically, the pairs of entering strands (at left) alternate between the positive element being assigned to the lower strand and the negative element being assigned to the lower strand. Thus, we will proceed by induction on four entering strands at a time.  We can directly verify by figure~\ref{figacol} that the lemma holds for $THK(3,2)$ and $THK(5,2)$
with $S_1 = 1$, $S_2 = 3$, $S_3 = 8$, and $S_4 = 20$.  Assume that, the lemma holds for $THK(4k-5 ,2)$ and $THK(4k-3,2)$ where $k \geq 2$. That is, as in
figure~\ref{figkind} assume that the $THK(m,2)$ diagram has been colored with $S_{4k-7} = w$, $S_{4k-6} = x$, $S_{4k-5} = y$ and $S_{4k-4} = z$.
We then assign the next four strands as in the figure, so that
$S_{4k-3} = w' = 2z + y + 1$, $S_{4k-2} = x' = 2w' + z + 1$.
$S_{4k-1} = y' = 2x' + w' + 1$, and $S_{4k} = 2y' +x' +1$.
This shows that the lemma also holds for
$THK(4k-1,2)$ and $THK(4k+1,2)$. \bx

\medskip

Our goal is to show that every non-trivial coloring of $THK(m,2)$ assigns different colors to different strands. We begin by observing that the coloring with $a = 1$ has this property.

\begin{lemma}
\label{lema1}%

Let $m \geq 3$ be odd and let $P_m$, the $m$th Pell number, be the determinant of the Turk's Head Knot $THK(m,2)$.
The $THK(m,2)$ diagram admits a $P_m$-coloring 
that assigns different colors to different strands.

\end{lemma}

\pr We may color our braid as we have done starting with zero, but let us replace $a$ with $1$.  As in lemma~\ref{lemSn}, we will look at the ascending sequence $S_n$ obtained by taking the positive element of the pairs in figure 4.  Observe that, in the $THK(m,2)$ diagram, the largest color in absolute value is the $m$th entering strand.  Although there are vertical strands with higher absolute value, these are forced to be equivalent to existing strands of lesser absolute value.  Thus, the highest absolute value for a color in $THK(m,2)$ is $S_{m-2}$.  Observe that for $m=3$, $S_{1}=1< \frac32=\frac{P_{3}}{2}-1$  Again, for $m=5$, $S_{3}=8<\frac{P_{5}}{2}-1$.  Now assume that, for all $3 \leq i<k$, $S_{i-2}<\frac{P_{i}}{2}-1$.  Then,

\begin{center}
$S_{k-2}=2S_{k-3}+S_{k-4}+1<2(\frac{P_{k-1}}{2}-1)+(\frac{P_{k-2}}{2}-1)+1<
\frac{2P_{k-1}+P_{k-2}}{2}-1=\frac{P_{k}}{2}-1$
\end{center}

So the absolute value of the color of every strand in our coloring of $THK(m,2)$ is less than half the determinant of the knot. This means that in this particular coloring, no two strands are assigned the same color modulo $P_m$.
\bx

\medskip

Finally, to prove the theorem, it remains only to verify that every non-trivial $P_m$-coloring assigns different colors to different strands. 

\medskip

\pr (of theorem~\ref{thmKHm2})
Let $m \geq 3$ be odd. As shown in theorem~\ref{thmdetm2}, the
determinant of $THK(m,2)$ is the $m$th Pell number $P_m$.
If $P_m$ is not prime, then we are done.  So, we may assume $P_m$ is prime. 

As above, up to a parameter $a$, we have essentially
one way of coloring the diagram mod $P_m$.  By lemma~\ref{lema1}, when $a=1$, that coloring assigns different colors to 
different strands. More generally, we've shown that
for any value of $a$ the strands will be labeled with different multiples of $a$. So, suppose two different strands
are labeled with the colors $c_1a$ and $c_2a$. This means that in
the $a=1$ coloring, these two strands are labeled $c_1$ and $c_2$, and,
by lemma~\ref{lema1}, $c_1 \not\equiv c_2$. Then, if $1 \leq a < P_m$, it
follows that $c_1a \not\equiv c_2a$. So, for any choice of $a$ (beside $a  = 0$ which leads to a trivial coloring) and, therefore, for any 
$P_m$-coloring
of the diagram, different strands are assigned different colors. \bx

\section{Proof of Theorem~\ref{thmdet}}

\setcounter{thm}{3}

In this section we prove

\begin{thm}
Let $m \geq 3$ and $n \geq 3$ be relatively prime integers.
If $m = 3$ or $m$ is even, then the Turk's Head Knot $THK(m,n)$ 
has composite determinant.
\end{thm}

\setcounter{thm}{8}

\begin{figure}[ht]
\begin{center}
\includegraphics[width=4in]{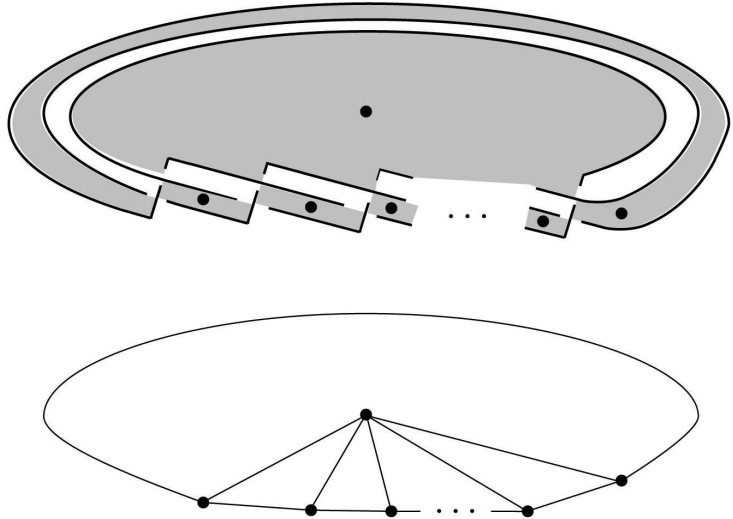}
\caption{The checkerboard graph of $THK(3,n)$ is a wheel on $n+1$ vertices.}\label{figwn}
\end{center}
\end{figure} 

\pr
If $m = 3$, the checkerboard graph of $THK(m,n)$ is a wheel of 
$n+1$ vertices (see figure~\ref{figwn}). The number of spanning trees, and hence the determinant of $THK(m,n)$ 
is $L_{2n} -2$ where $L_n$ is the $n$th Lucas number (see \cite{R}).
By induction, the determinant is $5 F_n^2$, where $F_n$ is the $n$th Fibonacci number, when $n$ is 
even and $L_n^2$ when $n$ is odd. Thus, the determinant is composite 
when $m = 3$.

\begin{figure}[ht]
\begin{center}
\includegraphics[width=4.5in]{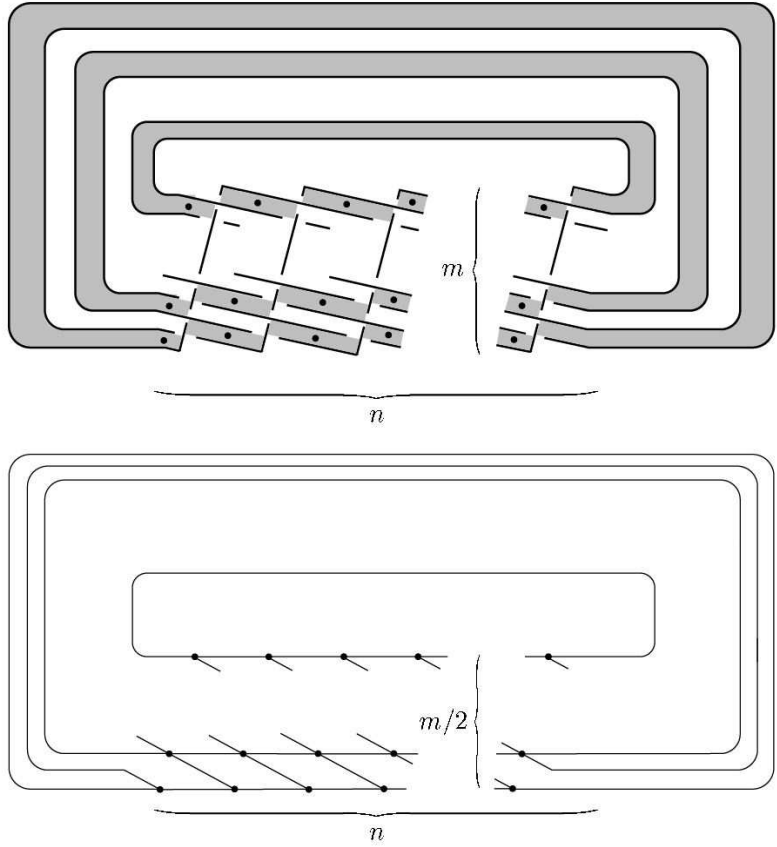}
\caption{For $m$ even, the checkerboard graph of $THK(m,n)$ is the tensor product of an $n$-cycle and a $m/2$-chain.}\label{figspanmn}
\end{center}
\end{figure} 

If $m \geq 4$ is even, the calculation of the number of spanning trees is
given in \cite{K}. In this case, (see figure~\ref{figspanmn}) the checkerboard graph is the tensor product of a cycle of $n$ vertices and a chain of $m/2$ vertices.  The number of spanning trees and hence the determinant of
$THK(m,n)$ is
$$S_{m,n} = n \prod_{\begin{array}{c}
1 \leq h \leq n-1 \\ 1 \leq k \leq \frac{m}{2} -1 \end{array}}
(4 \sin^2 \frac{h \pi}{n} + 4 \sin^2 \frac{k \pi}{m})$$
Now, the terms $4 \sin^2 \frac{h \pi}{n} + 4 \sin^2 \frac{k \pi}{m}$ constitute a complete set of conjugate roots of a polynomial with integral coefficients. 
So, their product is an integer. Moreover,  the product is greater than one. Thus, 
$S_{m,n}$ is of the form $n$ times an integer greater than one, and therefore is composite when $m$ is even. \bx

\section{A formula for the determinant $G_{m,n}$ of $THK(m,n)$}

\setcounter{thm}{4}

In this section we 
derive an expression $G_{m,n}$ that, we believe, is the determinant of $THK(m,n)$ when $m$ is odd. We then prove

\begin{thm}
$G_{m,n}$ is composite when $m,n \geq 3$ and $m$ is odd. 
\end{thm}

\setcounter{thm}{8}

Thus, if $G_{m,n}$ is in fact the determinant of these knots,
this would complete an argument that the Harary-Kauffman conjecture holds
for all Turk's Head Knots. 

Our formulation of $G_{m,n}$ is based on a 
matrix $\A_m$ we used to analyze coloring the $THK(m,n)$ braid. For example, if $m = 3$, the matrix is
$$
\A_3= \left( \begin{array}{ccc}
    2&-1&0\\
    0&0&1\\
    -1&0&2
\end{array} \right). $$

\begin{figure}[ht]
\begin{center}
\includegraphics[width=1.7in]{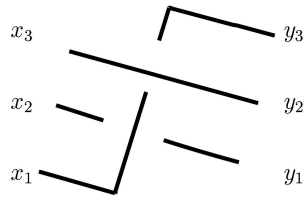}
\caption{The braid for $\A_3$.}\label{figA3}
\end{center}
\end{figure} 

The idea is that $THK(3,n)$ is
obtained by $n$ repetitions of the braid of figure~\ref{figA3}.
If we label the 
strands as shown, then, using the coloring relation $2x -y - z \equiv 0$,
we can view the strands $\vec{y} = (y_1, y_2, y_3)$ leaving at right as being derived from those entering at left,
$\vec{x} = (x_1,x_2,x_3)$, by matrix multiplication: $\vec{y} = \A_3 \vec{x}$. For the diagram of $THK(3,n)$ we repeat this multiplication $n$ times so that the vector leaving the braid at right is $\A_3^n \vec{x}$. 

However, in order to have a valid coloring we require the
numbers entering the braid to be the same as the numbers exiting the
braid as they are in fact the same strands. That is, we require that
$\A_m^n\vec{x}=\vec{x}$. In other words, $\vec{x}$ represents a valid
coloring of $THK(m,n)$ if it is an eigenvector of
$\A_m^n$ with eigenvalue 1. 

Now, any constant vector $\vec{x} = (c,c,c)$ will be an eigenvector of
eigenvalue one and corresponds to the trivial coloring where all
strands are assigned the same color $c$. So, the matrix $\A_n^m$ will 
always have 1
as an eigenvalue. If $THK(m,n)$ admits a (non-trivial) $p$-coloring, then, modulo $p$, $\A_n^m$ will have additional eigenvectors beyond the constant vectors. This means that 1 occurs as a root more than once in the characteristic
polynomial of $\A_n^m$ taken modulo $p$. Thus, we can discover the 
valid crossing moduli, and from these the knot's determinant, by looking at
the characteristic polynomial of $\A_n^m$. 

Below, we carry out this program to find a number $G_{m,n}$ such that 
the characteristic polynomial of $\A_n^m$ does, indeed, have 1
as a multiple root modulo $G_{m,n}$. However, this is only a necessary condition for
there to be a coloring mod $G_{m,n}$. In other words, we can show that the algebraic 
multiplicity of $1$ as an $\A_n^m$ eigenvalue is at least two when working modulo 
$G_{m,n}$. However, we can't be sure that the geometric multiplicity is
also greater than one and so we don't know for sure that there is valid coloring modulo $G_{m,n}$. 

On the other hand, computer experiments show that our formula for $G_{m,n}$ agrees with the determinant at least for $m \leq 23$ and $n \leq 29$.

Let us then proceed with the calculation of $G_{m,n}$.
For $m=2$, $3$, $4$, $5$, the matrices $\A_m$ have the following form:

\begin{displaymath}
\A_2 = \left( \begin{array}{cc}
    2&-1\\
    1&0
\end{array} \right)
\end{displaymath}\\

\begin{displaymath}
\A_3 = \left( \begin{array}{ccc}
    2&-1&0\\
    0&0&1\\
    -1&0&2
\end{array} \right)
\end{displaymath}\\

\begin{displaymath}
\A_4 = \left( \begin{array}{cccc}
    2&-1&0&0\\
    0&0&1&0\\
    -2&0&4&-1\\
    -1&0&2&0
\end{array} \right)
\end{displaymath}\\

\begin{displaymath}
\A_5 = \left( \begin{array}{ccccc}
    2&-1&0&0&0\\
    0&0&1&0&0\\
    -2&0&4&-1&0\\
    0&0&0&0&1\\
    1&0&-2&0&2
\end{array} \right)
\end{displaymath}\\

In general, for $m \geq 2$, $\A_{2m}$ has the form

\begin{tiny}

\begin{equation}\label{eqA2m}
\left( \begin{array}{cccccccccccc}
    2&-1&0&0&0&\cdots&\cdots&\cdots&0&0&0&0\\
    0&0&1&0&0&\cdots&\cdots&\cdots&0&0&0&0\\
    -2&0&4&-1&0&\cdots&\cdots&\cdots&0&0&0&0\\
    0&0&0&0&1&\ddots&\cdots&\cdots&0&0&0&0\\
    2&0&-4&0&4&\ddots&\ddots&\cdots&\vdots&\vdots&\vdots&\vdots\\
    \vdots&0&0&0&0&\ddots&\ddots&\ddots&\vdots&\vdots&\vdots&\vdots\\
    \vdots&\vdots&4&0&-4&\ddots&\ddots&\ddots&0&\vdots&\vdots&\vdots\\
    \vdots&\vdots&\vdots&0&0&\ddots&\ddots&\ddots&1&0&\vdots&\vdots\\
    \vdots&\vdots&\vdots&\vdots&4&\ddots&\ddots&\ddots&4&-1&0&\vdots\\
    \vdots&\vdots&\vdots&\vdots&\vdots&\ddots&\ddots&\ddots&0&0&1&0\\
    2(-1)^{m+1}&0&4(-1)^{m}&0&4(-1)^{m+1}&\cdots&\ddots&\ddots&-4&0&4&-1\\
    (-1)^{m+1}&0&2(-1)^{m}&0&2(-1)^{m+1}&\cdots&\cdots&\ddots&-2&0&2&0
\end{array} \right)
\end{equation}\\
\end{tiny}

while $\A_{2m+1}$ will be

\begin{tiny}

\begin{equation}\label{eqA2m+1}
\left( \begin{array}{ccccccccccccc}
    2&-1&0&0&0&\cdots&\cdots&\cdots&0&0&0&0&0\\
    0&0&1&0&0&\cdots&\cdots&\cdots&0&0&0&0&0\\
    -2&0&4&-1&0&\cdots&\cdots&\cdots&0&0&0&0&0\\
    0&0&0&0&1&\ddots&\cdots&\cdots&0&0&0&0&0\\
    2&0&-4&0&4&\ddots&\ddots&\cdots&\vdots&\vdots&\vdots&\vdots&\vdots\\
    \vdots&0&0&0&0&\ddots&\ddots&\ddots&\vdots&\vdots&\vdots&\vdots&\vdots
\\
   \vdots&\vdots&4&0&-4&\ddots&\ddots&\ddots&0&\vdots&\vdots&\vdots&\vdots
\\
   \vdots&\vdots&\vdots&0&0&\ddots&\ddots&\ddots&1&0&\vdots&\vdots&\vdots\\
    \vdots&\vdots&\vdots&\vdots&4&\ddots&\ddots&\ddots&4&-1&0&\vdots&\vdots
\\
    \vdots&\vdots&\vdots&\vdots&\vdots&\ddots&\ddots&\ddots&0&0&1&0&\vdots
\\
    2(-1)^{m+1}&0&4(-1)^{m}&0&4(-1)^{m+1}&\cdots&\ddots&\ddots&-4&0&4&-1&0\\
    0&\vdots&0&\vdots&0&\cdots&\cdots&\ddots&0&0&0&0&1\\
    (-1)^{m}&0&2(-1)^{m+1}&0&2(-1)^{m}&\cdots&\cdots&\cdots&2&0&-2&0&2
\end{array} \right)
\end{equation}\\
\end{tiny}

Let $g_m(x)$ denote the characteristic polynomial of $\A_m$. 
By direct calculation,
\begin{eqnarray*}
g_2 &=& 1 -2x + x^2\\
g_3 &=& 1-4x + 4x^2 -x^3\\
g_4 &=& 1-6x + 10x^2 - 6x^3 + x^4\\
g_5 &=& 1 -8x +20x^2 - 20x^3 + 8x^4-x^5
\end{eqnarray*}
Using the form of $\A_{2m}$ and $\A_{2m+1}$ given in equations~\ref{eqA2m}
and \ref{eqA2m+1}, one can show by induction that the $g_m(x)$ satisfy the recursion
\begin{eqnarray*}
g_{m+1} & = & (1-x) g_m - x g_{m-1} \\
& = & g_m - x(g_m + g_{m-1}) 
\end{eqnarray*}

As we have mentioned earlier, $1$ is always a root of these polynomials. 
It will be convenient to instead work with the polynomials 
$d_{m} = g_m/(1-x)$ which satisfy the same recursion. Then,
\begin{eqnarray*}
d_2 &=& 1 -x\\
d_3 &=& 1-3x + x^2\\
d_4 &=& 1-5x + 5x^2 - x^3 \\
d_5 &=& 1 -7x +13x^2 - 7x^3 + x^4
\end{eqnarray*}

This sequence of polynomials is closely related to the Delannoy numbers:
\begin{equation}\label{eqD}
\begin{array}{cccccccccccccc}
   & &  &  &  &  & 1 &  &  &  &  &  &  &  \\
   &  &  &  &  & 1 &  & 1 &  &  &  &  &  &  \\
   &  & &  & 1 &  & 3 &  & 1 &  &  &  &  &  \\
   &  &  & 1 &  & 5 &  & 5 &  & 1 &  &  &  &  \\
   &  & 1 &  & 7 &  & 13 &  & 7 &  & 1 &  &  &  \\
   & 1 &  & 9 &  & 25 &  & 25 &  & 9 &  & 1 &  &  \\
  1 &  & 11 &  & 41 &  & 63 &  & 41 &  & 11 &  & 1 &  \\
&&&&&&\vdots&&&&&&&
\end{array}
\end{equation}
The Delannoy numbers are defined by the recurrence
$D_{i,k} = D_{i-1,k} + D_{i-1,k-1} + D_{i-1,k-2}$. For example the $13$ in the fifth row of equation~\ref{eqD} is the sum of the three terms above it, namely $5$, $5$, and $3$. 

Thus, as may be verified by induction, the coefficients of the
polynomials $d_m(x)$ arising from the characteristic polynomials
of the $\A_m$ matrices are the Delannoy numbers:
$$d_m(x) = \sum_{k=0}^{m-1} (-1)^k D_{m-1,k} x^k.$$
It follows that the $d_m(x)$ is
palindromic in that $|D_{m-1,k}| = |D_{m-1,m-1-k}|$. Therefore, 
if $\alpha$ is a root of $d_{m}(x)$, then so too is $\alpha^{-1}$.
Moreover, we can argue that these roots are positive real numbers.

\begin{lemma} \label{lemapos}
For $m \geq 3$, the roots of $d_m$ are real and positive.  
\end{lemma}

\pr
In fact, we will argue that the roots of $g_m$ are real and positive. 
Note that,
since $|d_m(1)|=1$ when $m$ is odd, $1$ is a simple root of
$g_m$. Our argument will also show that $1$ is a double root
of $g_m$ when $m$ is even. Since $|d_m(1)|=0$ in this case,
we already know that $1$ is, at least, a double root.

We first analyze the properties of $g_3$ and $g_4$ and show by the
recursion relation that $g_5$ and $g_6$ must have the same
properties
(i.e., positive real roots and simple or double root at $1$); 
we then apply the same analysis to $g_{m-1}$ and $g_m$.

Let $s_i$ be the root of $g_i$ of minimal magnitude. The polynomial
$g_3$ has $s_3$ = $\frac{1}{\phi^2} \approx 0.381966$, and $g_4$ has
$s_4$ = $2-\sqrt 3 \approx 0.267949$.  (Here $\phi$ denotes the 
Golden Ratio.)
Note that, $s_4 < s_3$. Now
we proceed to $g_5$ and $g_6$. From our recursion relation, we have;
\begin{equation*}
g_5 = (1-x)g_4 -xg_3
\end{equation*}

Set $\tilde g_3$ = $x\frac{g_3}{1-x}$.  If $r \ne 1$ is a root of
$g_5$, then
\begin{equation}\label{real.pos.root}
g_4(r) = \tilde g_3 (r).
\end{equation}

\begin{figure}[ht]
\begin{center}
\includegraphics[width=2.8in]{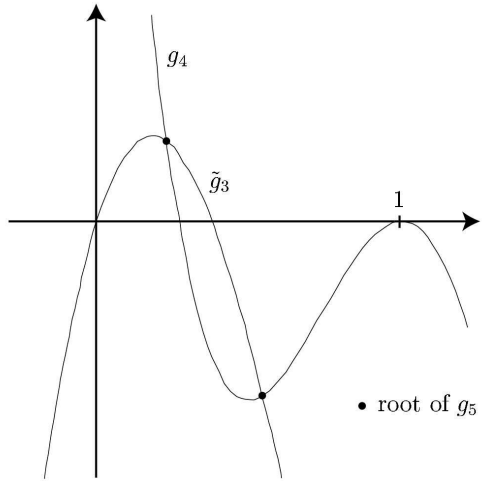}
\caption{Roots of $g_5$} \label{figg5}
\end{center}
\end{figure}

Since $x=1$ is a simple root of $g_3$, the polynomial $\tilde g_3$
has no root at $x=1$.  Also, $\tilde g_3(s_3) = \tilde g_3(0) = 0$.
We have $\tilde g_3(x)>0$ for $x \in (0, s_3)$, and since $s_4 \in
(0,s_3)$ and $g_4(0) > 0$, there must be some $r_1 \in (0,s_4)$ that
satisfies equation \eqref{real.pos.root}.  Also, $g_4(s_3) < 0 $ and
$\tilde g_3(x) < 0$ for $x \in (s_3, \frac{1}{s_3})$. Since $g_4(1) = 0$
and $\frac{1}{s_3} > 1$, there must be some $r_2 \in (s_3, 1)$
that satisfies equation \eqref{real.pos.root}; see figure $12$. The
points $r_1$, $\frac{1}{r_1}$, $r_2$, $\frac{1}{r_2}$ are roots of
$g_5$. Since $x=1$ is also a root, we have accounted for all the
roots of $g_5$ and they are all real and positive.

Consider $g_6$ = $(1-x)g_5 - xg_4$.  In $[0,1]$, every root of $g_4$
not equal to $1$ lies between two roots of $g_5$, so we can apply
the same argument to show that $\tilde g_4 $ = $x \frac{g_4}{1-x}$
intersects $g_5$ twice in the interval $(0,1)$. We conclude that
$g_6$ has two roots between $0$ and $1$, two roots at $1$, and two
roots after $1$.

Now assume the lemma holds for $m$ even and that $g_m$ and $g_{m-1}$
have only positive real roots. Also assume that these polynomials
have the property that in $[0,1]$, every root of $g_{m-1}$ lies
between two roots of $g_m$, and $s_m < s_{m-1}$.  Set $\tilde
g_{m-1} = x\frac{g_{m-1}}{1-x}$. Using the same analysis as above,
we know all the roots $r$ of $g_{m+1}$ not equal to $1$ satisfy $g_m(r)$ =
$\tilde g_{m-1}(r)$. The polynomial $g_m$ has $\frac{m-2}{2}$ roots
before $1$; again, using the same analysis as above we find that
$g_m$ and $\tilde g_{m-1}$ intersect $\frac{m-2}{2} + 1$ =
$\frac{m}{2}$ times, and these intersections occur between the roots
of $g_m$. Thus, $g_{m+1}$ has $\frac{m}{2}$ roots between $0$ and
$1$, $1$ root at $1$, and another $\frac{m}{2}$ roots after $1$
accounting for all its $m+1$ roots. Also $s_{m+1} < s_m$, and
continuing with the same analysis shows $g_{m+1}$ and $\tilde g_m$
intersect $\frac{m}{2}$ times between $0$ and $1$.  These
intersections account for $m$ of the $m+2$ roots of $g_{m+2}$, and
since $g_{m+2}$ has at least a double root at $x = 1$, it must be
that $g_{m+2}$ has exactly a double root at $x = 1$, and $m+2$ real
positive roots. \bx

\medskip

For the remainder of this section, we will assume $m$ is odd.
The roots of $d_m(x)$ are then the
pairs $\alpha_i$, $\alpha_i^{-1}$, $i = 1, \ldots ,
(m-1)/2$ and we define
$$G_{m,n} = \left[
\prod_{i=1}^{(m-1)/2} (\alpha_i^{n/2} - \alpha_i^{-n/2}) \right]^2$$

Let $d_m^{(n)}(x)$ denote the polynomial whose roots are the
$n$th powers of the roots of $d_m(x)$. In other words, 
the characteristic polynomial of $\A_m^n$ is $(1-x)d_m^{(n)}$.
Note that $G_{m,n} = d_m^{(n)}(1)$ so that, as promised,
$1$ is a multiple root of the characteristic polynomial of $\A_m^n$ modulo 
$G_{m,n}$.

It is also easy to verify that $G_{m,n}$ is the determinant of 
$THK(m,n)$ when $m=3$. Moreover, 
as in the case $m=3$, the $G_{m,n}$  are squares for
$n$ odd and of the form $G_{m,2}$ times a square when $n$ is even. 
This observation will allow us to show that $G_{m,n}$ is composite.

\medskip

\pr (of theorem~\ref{thmGmn}) 
First, let $n \geq 3$ be odd. In this case we have
\begin{eqnarray*}
G_{m,n} & = & \prod_i (\sqrt{\alpha_i} - \frac{1}{\sqrt{\alpha_i}})^2 \left[
\prod_i (\alpha_i^k + \alpha_i^{k-1} + \cdots + \alpha_i^{-k+1} + \alpha_i^{-k}) \right]^2 \\
& = & G_{m,1} (G'_{m,n})^2
\end{eqnarray*}
where $k = (n-1)/2$ and 
$$G'_{m,n} = \prod_i (\alpha_i^k + \alpha_i^{k-1} + \cdots + \alpha_i^{-k+1} + \alpha_i^{-k}).$$
Now, $|G_{m,1}| = |d_m(1)| = 1$ and $G'_{m,n}$ is also an integer as it
is a symmetric function of the roots $\alpha_i$. 
Also, for each $i$, $(\alpha_i^k + \alpha_i^{k-1} + \cdots + \alpha_i^{-k+1} + \alpha_i^{-k})$ is a positive real number greater than $1$ since, by lemma~\ref{lemapos},
$\alpha_i$ is a postive real number
and the sum includes $\alpha^0 = 1$. Therefore, 
$G_{m,n} = (G'_{m,n})^2$ is composite when $n \geq 3$ is odd.

Now, let $n \geq 4$ be odd. Again, we can factor
\begin{eqnarray*}
G_{m,n} & = & \prod_i ({\alpha_i} - \frac{1}{\alpha_i})^2 \left[
\prod_i (\alpha_i^k + \alpha_i^{k-2} + \cdots + \alpha_i^{-k+2} + \alpha_i^{-k}) \right]^2 \\
& = & G_{m,2} (G'_{m,n})^2
\end{eqnarray*}
where $k = (n-2)/2$ and 
$$G'_{m,n} = \prod_i (\alpha_i^k + \alpha_i^{k-2} + \cdots + \alpha_i^{-k+2} + \alpha_i^{-k}).$$

As in the $n$ even case, $G'_{m,n}$ is an integer greater than 1 and, so,
$G_{m,n}$ is also composite in this case.
\bx

\end{document}